\documentclass[reqno]{amsart}
\usepackage{amscd}
\newtheorem{define}{Definition}[section]
\newtheorem{cor}[define]{Corollary}
\newtheorem{prop}[define]{Proposition}
\newtheorem{teo}{Theorem}[section]
\theoremstyle{definition}
\newtheorem{ejem}[define]{Example}
\newtheorem{lema}[teo]{Lemma}
\theoremstyle{remark}
\newtheorem{nota}[teo]{Remark}

\begin{document}

\newcommand{\re}[1]{(\ref{#1})}

\newcommand{\DeliusA}{\cite{DeliusA}}
\newcommand{\Delius}{\cite{Delius}}
\newcommand{\DeliusB}{\cite{DeliusB}}
\newcommand{\Demidov}{\cite{Demidov}}
\newcommand{\Gurevich}{\cite{Gurevich}}
\newcommand{\GurevichA}{\cite{GurevichA}}
\newcommand{\Humphreys}{\cite{Humphreys}}
\newcommand{\Lyubashenko}{\cite{Lyubashenko}}
\newcommand{\Majid}{\cite{Majid}}
\newcommand{\Ringel}{\cite{Ringel}}
\newcommand{\Scheunert}{\cite{Scheunert} }
\newcommand{\ScheunertA}{\cite{ScheunertA} }
\newcommand{\Sudbery}{\cite{Sudbery} }
\newcommand{\Yamane}{\cite{Yamane}}
\newcommand{\Wolfram}{\cite{Wolfram}}

\newcommand{\acc}{\underline{\quad}\cdot\underline{\quad}}

\title[PBW for generalized Lie color algebras]{A Poincar\'e-Birkhoff-Witt theorem
for generalized  Lie color algebras}

\author[C. Bautista]{C\'esar Bautista}
\address{Facultad de Ciencias de la Computaci\'on.
Benem\'erita Universidad Aut\'onoma de Puebla.
Edif. 135. 14 sur y Av.
San Claudio. Ciudad Universitaria.
Puebla, Pue. C.P. 72570. M\'exico.}
\email{bautista@solarium.cs.buap.mx}

\begin{abstract}
A proof of Poincar\'e-Birkhoff-Witt theorem is given for
a class of generalized Lie algebras  closely related to the Gurevich $S$-Lie algebras.
As concrete
examples, we construct the positive (negative) parts of the quantized
universal enveloping algebras of type $A_{n}$ and
$M_{p,q,\epsilon}(n,\mathbb{K})$, which is
a non-standard quantum deformation  of $GL(n)$. In particular,
we get, for both algebras, a unified proof of the
Poincar\'e-Birkhoff-Witt theorem and we show that they are genuine
universal enveloping algebras of certain generalized Lie algebras.
\end{abstract}

\maketitle

\bigskip

{\bf PACS numbers}:  02.10.Vr, 02.10.Tq, 02.20.Sv

{\bf Running title}: PBW for generalized Lie color algebras

\begin{section}{Introduction}
\begin{sloppypar}
In the paper \Yamane, H. Yamane presented a proof of the
Poincar\'e-Birkhoff-Witt (PBW) theorem
for some class of quantum groups:
Drinfeld-Jimbo quantum groups of type $A_{n}$. In his proof he did not
use explicitly the Lie algebra theory concepts.
\end{sloppypar}

In this paper we show that Yamane used in an implicit manner
some generalized Lie algebra. Such a generalized Lie algebra  will be called
{\it $T$-Lie algebra}.

The $T$-Lie algebras satisfy not only generalized antisymmetry and Jacobi
identity, but aditional properties like multiplicativity, (also generalized,
in the same way as the Gurevich $S$-Lie algebras \Gurevich).
Such $T$-Lie algebras
arise in a natural way embedded in the positive and negative parts of
the Drinfeld-Jimbo quantum groups $U_{q}(sl_{n+1})$ of type $A_{n}$.

Our $T$-Lie algebras share some properties with the
$S$-Lie algebras. But they are not equivalent, for example,
$T$-Lie algebras satisfy a weaker  multiplicativity condition. In
particular, there are some $T$-Lie algebras which are not $S$-Lie algebras.
However, classical Lie algebras \Humphreys, Lie superalgebras \ScheunertA and
Scheunert generalized Lie algebras  (Lie color algebras) \Scheunert  are  all
$T$-Lie algebras.

These $T$-Lie algebras are related to the problem of finding the
appropriate definition
of {\it a quantum Lie algebra}. There are already some
generalized Lie algebras   proposed to solve this problem:
Majid braided Lie algebras \Majid, Delius-Gould quantum Lie
algebras \Delius, new generalized Lie algebras of
Gurevich-Rubstov \GurevichA, generalized Lie algebras due to
Lyubashenko-Sudbery \Lyubashenko, among others. But the Delius-Gould definition
and the Gurevich-Rubstov also, depends on the associated universal
enveloping algebra. This is not the case for the $T$-Lie
algebras. Our axioms imply the properties of the universal enveloping
algebra. In particular we shall prove the PBW theorem.

The generalized Lie algebras axioms of Lyubashenko-Sudbery are not
enough in order to obtain a PBW theorem (see example
\ref{L-S}). While the main difference with the braided Lie algebras
of Majid is that the symmetry of our $T$-Lie algebras is not a braid
morphism. Only  a part of such symmetry is braided.

In particular, we get a $T$-Lie algebra $(sl_{n+1}^{\pm})_{q}$
which is a deformation  of the Lie subalgebra of upper (lower)
triangular matrices. Such generalized Lie algebra meets almost all the
requirements of {\it a quantum Lie algebra} in the sense of
Lyubashenko-Sudbery \Lyubashenko, (only fails the point 7; actually the
universal enveloping algebra of $(sl_{n+1}^{\pm})_{q}$
has no a Hopf algebra structure, but it seems
possible to define a {\it braided} Hopf algebra on it, however we do not try such
matter in this paper). Moreover, the universal enveloping of
$(sl_{n+1}^{\pm})_{q}$ is $U_{q}^{\pm}(sl_{n+1})$ the positive part of the
Drinfeld-Jimbo quantum group of type $A_{n}$, therefore the diagram in
Figure \ref{cd1} commutes.
\bigskip
\begin{figure}[h]
\[\begin{CD}
sl_{n+1}^{\pm}      @>>> (sl_{n+1}^{\pm})_{q} \\
     @VVV                     @VVV \\
  U(sl_{n+1}^{\pm}) @>>> U_{q}^{\pm}(sl_{n+1})=U\,(sl_{n+1}^{\pm})_{q}
\end{CD}
\]
\caption{The classical $sl_{n+1}^{\pm}$ and the quantum $(sl_{n+1}^{\pm})_{q}$}
\label{cd1}
\end{figure}
This means that, relative to $U_{q}^{\pm}(sl_{n+1})$, the $T$-Lie algebra
$(sl_{n+1}^{\pm})_{q}$  satisfies, in some sense, the quantum Lie algebra
condition of Delius \DeliusA.

Some possible physical applications of the formalism of generalized Lie
algebras are in the affine Toda theories \DeliusB, quantum integrable
systems \DeliusB, and gauge theory \Lyubashenko.

The paper is organized as follows. In Sec. \ref{def} we shall define the
$T$-Lie algebras. In Sec. \ref{ejem} a list of classical and  new Lie
algebras is given. In Sec. \ref{UEA} we shall define the universal
enveloping algebra of a $T$-Lie algebra  and we shall prove that expecting an
analogue at PBW theorem for any such universal enveloping algebra constructed
by means of commutators is too much, we have to restrict our generalized Lie
algebras in an adecuate way. However, in Sec. \ref{rel} we persuit
the classical idea  to prove the PBW theorem \Humphreys by constructing  a representation of the
universal enveloping algebra on
the symmetric algebra (with modifications inspirated by \Yamane).
In Sec. \ref{rep} the definition of a representation of $T$-Lie algebra
is given. In Sec. \ref{PBW} we shall prove
the PBW theorem for the universal enveloping of an
adequate $T$-Lie algebra. Some remarks about braid morphisms are
given in Sec. \ref{braids}. The Sec. \ref{nonstandard} is devoted to explain why we can apply
the $T$-Lie algebras theory to a non-standard quantum deformation
algebra \Demidov of $GL(n)$.  Similar explanations are given in Sec.
\ref{A+} but now dealing with $U_{q}^{\pm}(sl_{n+1})$ the positive (negative)
parts of the Drinfeld-Jimbo quantum groups of type $A_{n}$. In
particular, in Sec. \ref{A+} we shall prove that $U_{q}^{\pm}(sl_{n+1})$ is a
genuine universal enveloping algebra of certain $T$-Lie algebra.
\end{section}

\begin{section}{The notion of $T$-Lie algebra}\label{def}
Let $k$ be a commutative unitary ring.

\begin{define}
\begin{sloppypar}
A $k$-algebra $A$ is {\it strictly graded} if there exist
$k$-submodules $(A_{\eta})_{\eta\in\mathbb{N}}$
such that
\end{sloppypar}
\noindent
\[
A=\oplus_{\eta\in\mathbb{N}}A_{\eta}\text{ and }A_{\eta_{1}}\cdot
A_{\eta_{2}}\subseteq A_{\eta_{1}+\eta_{2}-1}
\]
for all $\eta_{1},\eta_{2}\in\mathbb{N}$. For $a\in
A_{\eta}$, we shall put $\eta(a)=\eta.$
\end{define}

\begin{nota}
For such graded algebras $A$ we can induce a filtration of
$A\otimes_{k}A$ given by
\[
(A\otimes_{k}A)_{\eta} =\oplus_{\eta_{1}+\eta_{2}\leq\eta}A_{\eta_{1}}\otimes
A_{\eta_{2}}
\]
\end{nota}

 Let $L$ be a free $k$-module with a given basis
$\mathcal{B}$ totally ordered.

\begin{define}
Denote by $L^{n}$ the $k$-submodule  of $L^{n\otimes}$
generated by
\[
x_{i_{1}}\otimes\ldots\otimes x_{i_{n}},\quad x_{i_{1}}<\ldots
<x_{i_{n}},\quad(x_{i_{j}}\in \mathcal{B}),
\]
and by $\,^{n}L$ the $k$-submodule generated by
\[
x_{i_{1}}\otimes\ldots\otimes x_{i_{n}},\quad x_{i_{1}}>\ldots
>x_{i_{n}},\quad (x_{i_{j}}\in \mathcal{B}).
\]
\end{define}

\begin{define}

The module $L$ with $k$-morphisms
\begin{align}
S:L\otimes_{k}L&\rightarrow L\otimes_{k} L, \text{\it
 (presymmetry)}\\
T:L\otimes_{k}L&\rightarrow L\otimes_{k} L, \text{\it
 (symmetry)}\\
\langle,\rangle: L\otimes_{k}L&\rightarrow
L\otimes_{k}L, \text{\it  (pseudobracket)}\\
[,]:L\otimes_{k}L&\rightarrow L, \text{\it  (bracket)}
\end{align}
is called {\it $T$-Lie algebra with basis
$\mathcal{B}$ (or basic $T$-Lie algebra)} if, for
$S_{12}=S\otimes_{k}Id_{L},\,S_{23}=Id_{L}\otimes_{k} S$, the
following axioms are satisfied:

\begin{enumerate}
 \item
  \begin{enumerate}
  \item $S^{2}=Id$
  \item $S(x\otimes y)=q_{x,y}y\otimes x$, for certain $q_{x,y}\in k$,
  $\forall\,x,y\in \mathcal{B}$
  \item  {\it (Multiplicativity)}
    \begin{enumerate}
    \item
    $S(Id\otimes_{k}[,])|_{L^{3}}
    =([,]\otimes_{k}Id)S_{23}S_{12}|_{L^{3}}$
    \item $S([,]\otimes_{k} Id)|_{L^{3}}
    =(Id\otimes_{k}[,])S_{12}S_{23}|_{L^{3}}$
    \end{enumerate}
  \end{enumerate}
 \item {\it (Stability)}
    \begin{enumerate}
    \item  \label{est1} There exists a strict grading
           \[
            L=\oplus_{\eta}L_{\eta}
           \]
           of $L$ relative to $[,]$.
    \item \label{est2}
    \[
    \langle L_{\eta_{1}}\otimes L_{\eta_{2}}\rangle\subseteq
    (\,L\otimes_{k}L\,)_{\eta_{1}+\eta_{2}-1}
    \]
    for all $L_{\eta_{1}},\,L_{\eta_{2}}$.
        \end{enumerate}
 \item $ T=S +  \langle,\rangle$
  \item {\it (Antisymmetry)}\begin{enumerate}
        \item $[,]T=-[,]$
        \item $\langle,\rangle S=-\langle,\rangle$
        \item  $[,]\langle,\rangle=0$
        \end{enumerate}
\item {\it (Jacobi Identity)}
\[
[,](\,(Id\otimes_{k}[,])S_{12}S_{23}-
([,]\otimes_{k}
Id)S_{23}S_{12}+(Id\otimes_{k}[,])S_{23}S_{12}\,)|_{^{3}L}=0
\]
\end{enumerate}
\end{define}
Multiplicativity conditions are to control commutation relations in the
universal enveloping algebra,
whereas stability conditions are to obtain a good
gradation in the corresponding symmetric algebra.
\begin{define}
Let $L_{i}$ be a basic $T$-Lie  algebra with bracket $[,]_{i}$, pseudobracket
$\langle_,\rangle_{i}$ and presymmetry $S_{i}$, $i=1,2$. A
$k$-morphism $f:L_{1}\rightarrow L_{2}$ is called {\it a $T$-Lie
morphism} if $f$ is a morphism of graded algebras relative to $[,]_{i}$,
$i=1,2$ and the
diagrams in the Figure \ref{morph} commute.
\begin{figure}[ht]
\[
\begin{CD}
L_{1}\otimes_{k}L_{1} @>{\langle,\rangle}_{1}>> L_{1}\otimes_{k}L_{1}  \\
@V{f\otimes f}VV                        @VV{f\otimes f}V \\
L_{2}\otimes_{k}L_{2} @>{\langle,\rangle_{2}}>> L_{2}\otimes_{k}L_{2}
\end{CD}
\quad
\begin{CD}
L_{1}\otimes_{k}L_{1} @>S_{1}>> L_{1}\otimes_{k}L_{1}  \\
@V{f\otimes f}VV                        @VV{f\otimes f}V \\
L_{2}\otimes_{k}L_{2} @>{S_{2}}>> L_{2}\otimes_{k}L_{2}
\end{CD}
\]
\caption{A $T$-Lie algebra morphism}\label{morph}
\end{figure}
\end{define}
\end{section}

\begin{section}{Examples}\label{ejem}
In order to obtain a graduation in the stability conditions it suffices
to define a map $\eta:\mathcal{B}\rightarrow \mathbb{N}$ having properties
\eqref{est1} and \eqref{est2} in the stability axiom.
This remark will be used in the following examples.

\begin{subsection}{Some common Lie algebras}
\begin{ejem}
Classical Lie algebras over fields are basic $T$-Lie  algebras:
\[
[,]\text{ classical bracket,}\:\langle,\rangle=0,\,T=S\text{ usual swicht },\,\eta=1.
\]

\end{ejem}

\begin{ejem}
Lie superalgebras over fields \ScheunertA are basic $T$-Lie algebras:

Let $L=L_{0}\oplus L_{1}$ be a Lie superalgebra with bracket $[,]$ over a
field $k$. Let $\mathcal{B}_{\alpha}$ basis of $L_{\alpha}$, $\alpha=0,1$.
Define $S:L\otimes_{k}L\rightarrow L\otimes_{k}L$  on the basis
$\mathcal{B}=\mathcal{B}_{0}\cup \mathcal{B}_{1}$, by $S(x\otimes
y)=(-1)^{\alpha\beta}y\otimes x$ if $x\in L_{\alpha},\,y\in L_{\beta}$.
Besides
\[
\langle,\rangle=0,\,T=S,\,\eta=1.
\]
\end{ejem}

\begin{ejem}[Lie color algebras]
Let $k$ be a field of characteristic zero. Let
\[
L=\oplus_{\gamma\in\Gamma}L_{\gamma}
\]
be a $\epsilon$ Lie algebra \Scheunert, where $\Gamma$ is an abelian
group and $\epsilon$ is a commutation factor on $\Gamma$. Let
$[,]$ be bracket of $L$. Put $\mathcal{B}_{\gamma}$ basis of
$L_{\gamma}$ for each $\gamma\in\Gamma$.

Define 
\[
S(x\otimes y)=\epsilon(\alpha,\beta)y\otimes x,\text{ if
}x\in\mathcal{B_{\alpha}},y\in\mathcal{B_{\beta}},
\]
besides $\langle,\rangle=0$ and $\eta=1$.

Multiplicativity conditions follow easily from the definition of
commutation factor. We conclude that every $\epsilon$ Lie algebra is a
$T$-Lie algebra.
\end{ejem}
\end{subsection}

\begin{subsection}{Linear $T$-Lie algebras}
\begin{ejem}
\begin{sloppypar}
Let $e_{ij},\,1\leq i,j\leq n$ be standard basis of $gl_{n}$ matrices
$n\times n$ over a field $\mathbb{K}$. Let $[,]$ be the usual bracket in $gl_{n}$,
$sl_{n}^{+}$ the Lie subalgebra of upper triangular  matrices having trace
zero. We put
$x_{i}=e_{i,i+1},\,i=1,\ldots,n-1,\,x_{n}=[x_{1},x_{2}],\,x_{n+1}=[x_{2},x_{3}],
\ldots,
x_{2n-3}=[x_{n-2},x_{n-1}],x_{2n-2}=[x_{1},x_{n+1}]\ldots,
x_{3n-6}=[x_{n-3},x_{2n-3}],x_{3n-5}=[x_{1},x_{2n-1}],
x_{3n-4}=[x_{2},x_{2n}],\ldots,x_{m}=[x_{1},x_{m-1}]$ where
$m=n(n-1)/2$, besides we define $h_{i}=[x_{i},x_{i}^{t}],\,i=1,\ldots,m$
 diagonal matrices  in $sl_{n}$. Further,
$q=exp(t)\in
\mathbb{K}[[t]]$ formal series ring with indeterminate $t$ and
coefficents in $\mathbb{K}$, $k=\mathbb{K}[q,q^{-1}]$, $c_{i,j}\in \mathbb{Z}$ such
that $[h_{i},x_{j}]=c_{i,j} x_{j},\,1\leq i,j\leq m$.
\end{sloppypar}
Let $(sl_{n}^{+})_{q}$ be a free $k$-module with basis $\mathcal{B}=\{x_{i}\,|\,
1\leq i\leq m\}$. We may define an order in $\mathcal{B}$ according
to the Figure \ref{diag},
\begin{figure}[b]
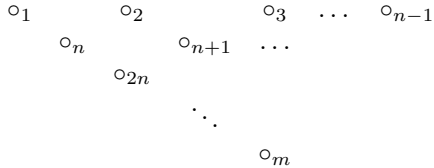

$$
\begin{array}{ccccccc}
\circ_{1} &             & \circ_{2} &  &\circ_{3} & \ldots & \circ_{n-1}\\
          & \circ_{n}   &           & \circ_{n+1} & \ldots \\
          &             & \circ_{2n}                        \\
          &             &           & \ddots      &         \\
          &             &           &             & \circ_{m}
\end{array} 
$$
\caption{The basic $T$-Lie algebra $(sl_{n+1}^{+})_{q}$}\label{diag}
\end{figure}
from left to right and up to bottom. For example
$x_{1}<x_{n}<x_{2}<x_{2n}$. The first time that a diagram (Auslander-Reiten quiver of type $A_{n-1}$)
of this type
appears related to quantum groups, is in Ringel's work about the
relationship between Poincar\'e-Birkhoff-Witt bases, quantum groups and
Hall algebras \Ringel.

Define:
\[
[x,y]_{q}=[x,y] \text{ if }x<y\in\mathcal{B}
\]
\[
\langle e_{ij},e_{uv}\rangle =\begin{cases}
                (q-q^{-1})e_{iv}\otimes e_{uj} & \text{if }i<u<j,u<j<v,\\
                 (q^{-1}-q)e_{uj}\otimes e_{iv} &\text{if
                 }u<i<v,i<v<j\\
                 0,\,&\text{otherwise.}
                \end{cases}
               \]
\[
S(x_{i}\otimes x_{j})=q^{c_{i,j}} x_{j}\otimes x_{i},\,\text{ if
}x_{i}<x_{j},
\]
\[
T=S+\langle,\rangle
\]
Finally,  we define $\eta$ in such way that every basic element in the
Figure \ref{diag} is in correspondence with a number belonging to the
Figure \ref{grad},
\begin{figure}[t]
$$
\begin{array}{ccccccc}
 1 &             & 2 &  & 3 & \ldots & n\\
   & 2           &   & 4&  & \ldots & 2(n-1)\\
   &             & 3 &  & 6 &\dots & 3(n-2)\\
   &             &   &\ddots& \\
   &             &   &      &n
\end{array}
$$
\caption{The graduation of $(sl_{n+1}^{+})_{q}$}\label{grad}
\end{figure}
this yields, $\eta(e_{ij})=i(j-i),\quad\forall i,j$.

The multiplicativity condition follows from properties
\[
x_{i}<x_{j}<x_{l}\Rightarrow\begin{cases}
                             [x_{i},x_{j}]<x_{l},\text{if
                             }[x_{i},x_{j}]\neq 0,\\
                             x_{i}<[x_{j},x_{l}],\text{if
                             }[x_{j},x_{l}]\neq 0,
                            \end{cases}
\]
and
\[
[h_{i},[x_{j},x_{l}]]=(c_{i,j}+c_{i,l})[x_{j},x_{l}]
\]

In the cases $n=2,3,4,5$, the Jacobi identity for $[,]_{q}$ can be verified by
straighforward calculations. We get that $(sl_{n}^{+})_{q}$ with bracket
$[,]_{q}$ is a basic $T$-Lie algebra, $n=2,3,4,5$.

Similarly, we can define
$(sl_{n}^{-})_{q}$.
\end{ejem}
\begin{nota}We get
\begin{equation*}
(sl_{n}^{\pm})_{q}|_{t=0}=sl_{n}^{\pm},
\end{equation*}
so in the cases $n=2,3,4,5,$ $(sl_{n}^{\pm})_{q}$  is a deformation of
$sl_{n}^{\pm}$ in the category of $T$-Lie algebras.  Later, in the
section \ref{A+}, such property will be generalized  for every $n$.
\end{nota}

\begin{ejem}
Starting from  $(sl_{4}^{+})_{q}$  we are going to build a new basic $T$-Lie
algebra, denoted $(\widetilde{sl_{4}^{+}})_{q}$. Its structure is:
\[
\widetilde{[,]}_{q}=[,]_{q},\quad\tilde{\langle,\rangle}=0,\quad
\tilde{S}=S,\quad\tilde{T}=S,\quad\tilde{\eta}=\eta.
\]
\end{ejem}

\begin{ejem}[ Non-standard quantum deformations \Demidov of $GL(n)$ ]\label{nonstand}
Let $\;p,\,q\;$ be units in a commutative unitary ring $k$ with $pq\neq 1$ and
choose $n(n-1)/2$ discrete parameters $\epsilon_{ij}$,
$\epsilon_{ij}=$$\pm 1,$ $1\leq i<j\leq
n,$ $\epsilon_{ii}=1,$ $\epsilon_{ji}=\epsilon_{ij}$.

The $k$-module $L_{p,q,\epsilon}(n,k)$ is then defined to be the free
$k$-module with basis
\[
\mathcal{B}=\{Z_{i}^{j}\,|\,1\leq i,j\leq n\}.
\]
We ordered $\mathcal{B}$ by putting $Z_{i}^{j}>Z_{u}^{v}$ if either
$i>u$, or $i=u$ and $j>v$. Define $\eta$ by
$\eta(Z_{j}^{i})=j\,3^{i-1}$.
Besides, we put $[,]=0$, and if $Z_{i}^{l}>Z_{u}^{v}$,
\begin{equation}
\langle Z_{i}^{l}, Z_{u}^{v}\rangle =\begin{cases}
                                     \epsilon_{vl}(p-q^{-1})Z_{i}^{v}\otimes
                                     Z_{u}^l,& \text{ if $i>u$ and $l>v$} \\
                                     0, &\text{ otherwise.}
                                     \end{cases}\label{pseudo}
\end{equation}
\newcommand{\cta}{Z_{u}^{v}\otimes Z_{i}^{l}}
\[
S(Z_{i}^{l}\otimes Z_{u}^{v})=\begin{cases}
                              \epsilon_{vl} p\cta,&\text{ if }i=u,l>v,\\
                              \epsilon_{ui}q\cta, &\text{ if }l=v,i>u,\\
                              \epsilon_{iu}\epsilon_{vl}p^{-1}q\cta,&\text
                              { if }i>u,v>l \\
                              \epsilon_{vl}\epsilon_{ui}\cta,&\text{ if
                              }i>u,l>v
                              \end{cases}
\]

To prove that $L_{p,q,\epsilon}(n,k)$ is a basic $T$-Lie algebra, since
\eqref{pseudo} it suffices to check the stability condition \eqref{est2} for
$Z_{i}^{l}>Z_{u}^{v}$, such that $i>u$ and $l>v$:
\[
(i-u)3^{l-1}>(i-u)3^{v-1}
\]
then
\[
i 3^{l-1}+u 3^{v-1}>i 3^{v-1}+u 3^{l-1}
\]
but this equation has left side
$\eta(Z_{i}^{l})+\eta(Z_{u}^{v})$ whereas the right side is
$\eta(Z_{i}^{v})+\eta(Z_{u}^{l})$. This proves stability conditions.
\end{ejem}
\end{subsection}
\end{section}
\begin{section}{Universal Enveloping Algebras}\label{UEA}
\begin{subsection}{Construction of $U(L)$}
\begin{define}
Let $L$  be a $T$-Lie algebra with basis $\mathcal{B}$, and
$\otimes_{k}L$ the $k$-tensor algebra of the module $L$.
The {\it universal enveloping algebra} $U(L)$ is the quotient
\[
U(L)=\otimes_{k}L/J
\]
where $J$ is the two sided ideal  generated by
\[
x\otimes y -T(x\otimes y)-[x,y],\: x,y\in \mathcal{B}
\]
\end{define}

Because the stability axiom \eqref{est2}, the algebra $U(L)$ have a
similar structure
to a {\it quadratic algebra with an ordering alghorithm \Sudbery.}
\end{subsection}
\begin{subsection}{Examples}
\begin{ejem}
\[
U(L_{p,q,\epsilon}(n,k))=M_{p,q,\epsilon}(n,k)
\]
is a non-standard quantum deformation \Demidov of $GL(n)$.
\end{ejem}
\begin{ejem}
\[
U(sl_{n+1}^{\pm})_{q}=U_{q}^{\pm}(sl_{n+1})
\]
positive (negative) part of the Drinfeld-Jimbo quantum group of type
$A_{n},\:n=3,4$.
\end{ejem}
\begin{ejem}\label{L-S}
\begin{sloppypar}
In $U(\widetilde{sl_{4}^{+}})_{q}$ the equation $x_{2}x_{6}=0$ holds.
Then $U(\widetilde{sl_{4}^{+}})_{q}$ is a enveloping algebra
where the  PBW theorem does not hold. So, if we want
a good enveloping algebra we have to add conditions to the $T$-Lie
algebras.
\end{sloppypar}

Besides, if $\beta$, $\tilde{S}$ denotes the bracket and the
symmetry of $\widetilde{(sl_{4}^{+}})_{q}$ respectivily and the
characteristic of the field $\mathbb{K}$ is zero, then for
$\gamma=Id-\tilde{S}$, where $Id$ is the identity morphism on
$\widetilde{(sl_{4}^{+}})_{q}^{2\otimes}$, the condition
$\gamma(t)=0$ implies $\beta(t)=0$. Moreover, if $\mathcal{B}$ is the
canonical basis of $\widetilde{(sl_{4}^{+}})_{q}$, and
$x,y,z$ are arbitrary elements in $\mathcal{B}$,  straighforward
calculations (using $Mathematica$ \Wolfram) gives:
\begin{eqnarray*}
\beta(\,\beta(x,y),z\,)&=
                   \beta(\,x,\beta(y,z)\,)-q_{x,y}\beta(\,y,\beta(x,z)\,)\\
\beta(\,z,\beta(x,y)\,)&=
                   \beta(\,\beta(z,x),y\,)-q_{x,y}\beta(\,\beta(z,y),x\,)
\end{eqnarray*}
where $\tilde{S}(x\otimes y)=q_{x,y}y\otimes x$. This means that
$\widetilde{(sl_{4}^{+}})_{q}$ has a structure of {\it balanced generalised
Lie algebra \Lyubashenko} and its universal enveloping algebra as such
generalised Lie algebra is the same as $T$-Lie algebra. Therefore, the
generalised Lie algebra axioms of Lyubashenko-Sudbery are not enough in
order to obtain a PBW theorem.
\end{ejem}
\begin{ejem}
Let $L$ be a basic $T$-Lie algebra. We are going to define a new $T$-Lie
algebra  $L^{0}$: $L^{0}=L$ in its structure of $k$-module,
$[,]^{0}=0,\langle,\rangle^{0}=0,\eta^{0}=\eta,S^{0}=S$ and
define
\[
\mathcal{S}(L)=U(L^{0})
\]
\begin{sloppypar}
$\mathcal{S}(L)$ is a free $k$-module with basis the monomials
formed
by the products
$z_{i_{1}}z_{i_{2}}\ldots z_{i_{r}}$ of elements of $\mathcal{B}$ such
that $r\geq 0, z_{i_{1}}\leq
z_{i_{2}}\leq\ldots \leq z_{i_{r}}\text{ where }z_{i_{j}}=z_{i_{j+1}}\text{ if
}q_{z_{i_{j}}z_{i_{j+1}}}=1$.
\end{sloppypar}
 Such an object
 $\mathcal{S}(L)$ will be called the {\it
$q$-symmetric algebra of $L$}.
\end{ejem}
\end{subsection}
\end{section}

\begin{section}{The relationship between universal enveloping algebras and
symmetric algebras}\label{rel}
\begin{sloppypar}
Let $L$  be a $T$-Lie algebra with basis $\mathcal{B}$,
$x_{\lambda}\in\mathcal{B}$,
$\Sigma=(x_{\lambda_{1}},\ldots,x_{\lambda_{u}})$ finite non-decreasing
sequence of elements of
$\mathcal{B}$. We write $z_{\lambda}=x_{\lambda}\in
\mathcal{S}(L)$, $z_{\Sigma}=z_{\lambda_{1}}\ldots z_{\lambda_{u}}\in
\mathcal{S}(L)$, $z_{\emptyset}=1\in\mathcal{S}(L)$,
$\eta(\lambda)=\eta(x_{\lambda}),\,\eta(\Sigma)=\eta(z_{\Sigma})=\eta(x_{\lambda_{1}})+\ldots
+\eta(x_{\lambda_{n}})$. Besides we put $x_{\lambda}\leq \Sigma$
if $x_{\lambda}\leq x_{\lambda_{1}}$.
\end{sloppypar}
\begin{lema}[A-B] Let $L$ be a $T$-Lie algebra with basis
$\mathcal{B}$. $\mathcal{P}=\mathcal{S}(L)$  $q$-symmetric algebra,
$\mathcal{P}_{p}$ the $k$-submodule generated by
$z_{\Sigma}$ such that $\eta(\Sigma)\leq p$.
There is a $k$-morphism
\begin{equation}
\acc:L\otimes_{k}\mathcal{P}\rightarrow \mathcal{P}
\end{equation}
satisfying
\begin{enumerate}
\item[(A)] $x_{\lambda}\cdot z_{\Sigma}=z_{\lambda}z_{\Sigma}$ for
$x_{\lambda}\leq \Sigma$;
\item[(B)] $x_{\lambda}\cdot z_{\Sigma}-z_{\lambda}z_{\Sigma}\in
\mathcal{P}_{\eta(\lambda)+\eta(\Sigma)-1}$.
\end{enumerate}
\end{lema}
\begin{proof}
By induction on $\eta(\lambda)+\eta(\Sigma)$. If
$\eta(\lambda)+\eta(\Sigma)=1$ then $\eta(\lambda)=1$ and $\Sigma=\emptyset$,
it follows $z_{\emptyset}=1$. Then define
\[
x_{\lambda}\cdot 1=z_{\lambda}
\]
so (A) and (B) holds.
Assume the existence of $x_{\lambda^{\prime}}\cdot z_{\Sigma^{\prime}}$
for
$\eta(\lambda^{\prime})+\eta(\Sigma^{\prime})<\eta(\lambda)+\eta(\Sigma)$
satisfying (A) and (B). We have to define  $x_{\lambda}\cdot
z_{\Sigma}$.

There are two cases: $\lambda\leq \Sigma$ or $\lambda\not\leq\Sigma$.

\underline{Case} $\lambda\leq \Sigma$: Because (A):
\[
x_{\lambda}\cdot z_{\Sigma}=z_{\lambda}z_{\Sigma}
\]

\underline{Case} $\lambda\not\leq \Sigma$: We may write
$\Sigma=(x_{\mu},N)$ with $x_{\mu}\leq N$ and $x_{\lambda}>x_{\mu}$.
Since
$\eta(N)<\eta(\Sigma)$ and because at the induction hypothesis
$x_{\lambda}\cdot z_{N}$
is already defined, and
\[
w=x_{\lambda}\cdot
z_{N}-z_{\lambda}z_{N}\in\mathcal{P}_{\eta(\lambda)+\eta(N)-1}.
\]
Moreover, from
$\eta(\mu)+\eta(\lambda)+\eta(N)-1<\eta(\mu)+\eta(\lambda)+
\eta(N)=\eta(\lambda)+\eta(\Sigma)$ it follows that
$x_{\mu}\cdot w$ is already defined.

We have
\begin{equation}
T(x_{\lambda}\otimes x_{\mu})=q_{\lambda\mu}x_{\mu}\otimes
x_{\lambda}+\sum_{i}\xi_{i}x_{\mu_{i}}\otimes x_{\lambda_{i}}
\label{sim}
\end{equation}
and  because at (B) and the induction hypothesis $x_{\lambda_{i}}\cdot
z_{N}\in\mathcal{P}_{\eta(\lambda_{i})+\eta(N)}$. As a consecuence
$x_{\mu_{i}}\cdot (\,x_{\lambda_{i}}\cdot z_{N}\,)$ is already
defined  because
$\eta(\mu_{i})+\eta(\lambda_{i})+\eta(N)<\eta(\lambda)+\eta(\mu)+\eta(N)$
according to stability axiom.

We may define
\begin{multline}
x_{\lambda}\cdot
z_{\Sigma}=q_{\lambda\mu}z_{\mu}z_{\lambda}z_{N}
+q_{\lambda\mu}x_{\mu}\cdot w
+\sum_{i}\xi_{i}x_{\mu_{i}}\cdot (x_{\lambda_{i}}\cdot z_{N})
+[x_{\lambda},x_{\mu}]\cdot z_{N}
\end{multline}
where $\;w=x_{\lambda}\cdot z_{N}-z_{\lambda}z_{N}$;
$\;[x_{\lambda},x_{\mu}]\cdot z_{N}$ is already defined because
$\eta([x_{\lambda},x_{\mu}])+\eta(N)<\eta(\lambda)+
\eta(\mu)+\eta(N)=\eta(\lambda)+\eta(\Sigma)$.

Now only remains to prove (B). From
$z_{\lambda}z_{\Sigma}=q_{\lambda\mu}z_{\mu}z_{\lambda}z_{N}$ we obtain
\[
x_{\lambda}\cdot
z_{\Sigma}-z_{\lambda}z_{\Sigma}=q_{\lambda\mu}x_{\mu}\cdot w
+\sum_{i}\xi_{i}x_{\mu_{i}}\cdot x_{\lambda_{i}}\cdot z_{N}
+[x_{\lambda},x_{\mu}]\cdot z_{N}.
\]
Moreover
\[
x_{\mu}\cdot w \in
\mathcal{P}_{\eta(\mu)+\eta(w)}=
\mathcal{P}_{\eta(\mu)+\eta(\lambda)+\eta(N)-1}=
\mathcal{P}_{\eta(\lambda)+\eta(\Sigma)-1}
\]
\[
x_{\mu_{i}}\cdot x_{\lambda_{i}}\cdot z_{N}\in
\mathcal{P}_{\eta(\mu_{i})+\eta(\lambda_{i})+\eta(N)}
\subseteq  \mathcal{P}_{\eta(\mu)+\eta(\lambda)-1+\eta(N)}
=\mathcal{P}_{\eta(\lambda)+\eta(\Sigma)-1}
\]
\[
[x_{\lambda},x_{\mu}]\cdot z_{N}\in
\mathcal{P}_{\eta(\lambda)+\eta(\mu)-1+\eta(N)}=
\mathcal{P}_{\eta(\lambda)+\eta(\Sigma)-1}
\]
imply
\[
x_{\lambda}\cdot z_{\Sigma}-z_{\lambda}z_{\Sigma}\in \
\mathcal{P}_{\eta(\lambda)+\eta(\Sigma)-1}
\]
\end{proof}

\begin{define}
Let $L$ be a $T$-Lie algebra with basis $\mathcal{B}$. We call
$L$  {\it adequate} if the morphism from the lemma (A-B)
is such that the condition
\begin{equation}
x_{\lambda^{\prime}}\cdot x_{\mu^{\prime}}\cdot z_{M}-
T(x_{\lambda^{\prime}}\otimes x_{\mu^{\prime}})\cdot z_{M}=
[x_{\lambda^{\prime}},x_{\mu^{\prime}}]\cdot z_{M}
\label{hipinducc}
\end{equation}
for all
$\eta(x_{\lambda^{\prime}})+\eta(x_{\mu^{\prime}})+\eta(M)\leq r$
implies
\begin{multline}
\langle x_{\lambda}, x_{\mu}\rangle\cdot x_{\gamma}\cdot
z_{N}-q_{\lambda\mu}[x_{\mu},[x_{\lambda},x_{\gamma}]]\cdot z_{N}=\\
q_{\mu\gamma}q_{\lambda\gamma}\langle x_{\gamma}[x_{\lambda},x_{\mu}]
\rangle \cdot z_{N}
+q_{\mu\gamma}q_{\lambda\gamma}x_{\gamma}\cdot
\langle x_{\lambda},x_{\mu}\rangle\cdot z_{N}\\
+q_{\mu\gamma}\langle x_{\lambda},x_{\gamma}\rangle\cdot x_{\mu}\cdot z_{N}
-q_{\lambda\mu}x_{\mu}\cdot\langle x_{\lambda},x_{\gamma}\rangle\cdot
z_{N} \\
+q_{\mu\gamma}[x_{\lambda},x_{\gamma}]\cdot x_{\mu}\cdot z_{N}
-q_{\lambda\mu}x_{\mu}\cdot [x_{\lambda},x_{\gamma}]\cdot z_{N}\\
+x_{\lambda}\cdot \langle x_{\mu},x_{\gamma}\rangle \cdot z_{N}
-q_{\lambda\mu}q_{\lambda\gamma}\langle x_{\mu},x_{\gamma}\rangle\cdot
x_{\lambda}\cdot z_{N}-q_{\lambda\gamma}q_{\lambda\mu}\langle
[x_{\mu},x_{\gamma}],x_{\lambda}\rangle\cdot z_{N}
\label{inducc}
\end{multline}
for every $x_{\lambda}>x_{\mu}>x_{\gamma}\in\mathcal{B}$, $x_{\gamma}\leq
z_{N}$ such that
$\eta(x_{\lambda})+\eta(x_{\mu})+\eta(x_{\gamma})+\eta(N)\leq r+1$.
\end{define}

\begin{lema}[C]
Let $L$ be an adequate $T$-Lie algebra with basis
$\mathcal{B}$, and $\mathcal{P}$ the related $q$-symmetric algebra.
Then there exists  a  $k$-morphism
$\underline{\quad}\cdot\underline{\quad}:L\otimes_{k}\mathcal{P}
\rightarrow\mathcal{P}$  such that
\begin{enumerate}
\item[(C)] $x_{\lambda}\cdot x_{\mu}\cdot z_{N}=T(x_{\lambda}\otimes x_{\mu})\cdot
z_{N}+[x_{\lambda},x_{\mu}]\cdot z_{N}$, $\forall
z_{N}\in\mathcal{P},\forall x_{\lambda},x_{\mu}\in\mathcal{B}$.
\end{enumerate}
\end{lema}

\begin{proof}
Let  $\acc$ be the morphism from lemma (A-B). There are two cases:
\begin{enumerate}
  \item\label{er} $\mu\leq N$ or $\lambda\leq N$;
  \item\label{do} $\mu\not\leq N$ and $\lambda\not\leq N$;
\end{enumerate}
\re{er}: Assume $\mu\leq N$ and $\mu<\lambda$. Let $M=(\mu,N)$, then,
by definition
\begin{align*}
x_{\lambda}\cdot x_{\mu}\cdot z_{N}&=x_{\lambda}z_{M}\text{ where
}\lambda\not\leq M \\
&=z_{\lambda}\cdot z_{M}+q_{\lambda\mu}x_{\mu}\cdot w +\langle
x_{\lambda},x_{\mu}\rangle \cdot z_{N}+[x_{\lambda},x_{\mu}]\cdot
z_{N}
\end{align*}
On the other hand,
\begin{align*}
T(x_{\lambda}\otimes & x_{\mu})\cdot z_{N}+[x_{\lambda},x_{\mu}]\cdot
z_{N}\\
&= q_{\lambda\mu}x_{\mu}\cdot x_{\lambda}\cdot z_{N}+\langle
x_{\lambda},x_{\mu}\rangle\cdot z_{N}+[x_{\lambda},x_{\mu}]\cdot z_{N}\\
&=q_{\lambda\mu}x_{\mu}\cdot
(z_{\lambda}z_{N})+q_{\lambda\mu}x_{\mu}\cdot w +\langle
x_{\lambda},x_{\mu}\rangle\cdot z_{N}+[x_{\lambda},x_{\mu}]\cdot z_{N}
\end{align*}
since $\mu<\lambda$ and $\mu\leq N$ it holds $z_{\lambda}z_{N}=c
z_{N^{\prime}}$  where $\mu\leq N^{\prime}$ and $c\in k$,
\[
x_{\mu}\cdot (z_{\lambda}z_{N})=cx_{\mu}\cdot
z_{N^{\prime}}=cz_{\mu}z_{N^{\prime}}=z_{\mu}z_{\lambda}z_{N},
\]
so
\[
q_{\lambda\mu}x_{\mu}\cdot(z_{\lambda}z_{N})=
q_{\lambda\mu}z_{\mu}z_{\lambda}z_{N}=z_{\lambda}z_{\mu}z_{N}=z_{\lambda}z_{M},
\]
therefore
\begin{align*}
T(x_{\lambda}\otimes x_{\mu})\cdot z_{N}&+[x_{\lambda},x_{\mu}]\cdot
z_{N}\\
&=z_{\lambda}z_{M}+q_{\lambda\mu}x_{\mu}\cdot w +\langle
x_{\lambda},x_{\mu}\rangle\cdot z_{N}+[x_{\lambda},x_{\mu}]\cdot
z_{N}\\
&=x_{\lambda}\cdot x_{\mu} \cdot z_{N}
\end{align*}
(i.e. (C) holds for $\mu<\lambda$).
It follows, multiplying by $-q_{\mu\lambda}$:
\begin{equation}
-q_{\mu\lambda}x_{\lambda}\cdot x_{\mu}\cdot
z_{N}=-x_{\mu}\otimes x_{\lambda}\cdot
z_{N}-q_{\mu\lambda}\langle x_{\lambda},x_{\mu}\rangle \cdot z_{N}
-q_{\mu\lambda}[x_{\lambda},x_{\mu}]\cdot z_{N}.
\label{mult}
\end{equation}
This implies, using antisymmetry,
\begin{equation}
x_{\mu}\cdot x_{\lambda}\cdot z_{N}-T(x_{\mu}\otimes x_{\lambda})\cdot
z_{N}=[x_{\mu},x_{\lambda}]\cdot z_{N}
\end{equation}
and we conclude that (C) also holds for $\lambda<\mu$.

\re{do}: Let $N=(\gamma,Q)$ where $\gamma\leq Q,\,
\gamma <\lambda,\,\gamma <\mu$. We proceed by induction on
$\eta(\lambda)+\eta(\mu)+\eta(N)$. Suppose that for each
$\eta(\lambda^{\prime})+\eta(\mu^{\prime})+\eta(N^{\prime})\leq r$ it holds (C). Then, for
$\eta(\lambda)+\eta(\mu)+\eta(N)\leq r+1$ we have:
\begin{equation}
x_{\mu}\cdot z_{N}=x_{\mu}\cdot(x_{\gamma}\cdot
z_{Q})=T(x_{\mu}\otimes x_{\gamma})\cdot
z_{Q}+[x_{\mu},x_{\gamma}]\cdot z_{Q}
\label{dblstar}
\end{equation}
because $\eta(\mu)+\eta(\gamma)+\eta(Q)=\eta(\mu)+\eta(N)\leq r$ and the
induction hypothesis.

Now, $x_{\mu}\cdot z_{Q}=z_{\mu}z_{Q}+w$ where
$w\in\mathcal{P}_{\eta(\mu)+\eta(Q)-1}$. We may apply (C) to
$x_{\lambda}\cdot x_{\gamma}\cdot(z_{\mu}z_{Q})$ since
$z_{\mu}z_{Q}=cz_{Q^{\prime}}$ where $c\in k$ and $\gamma\leq Q^{\prime}$
because $\gamma\leq Q$, $\gamma<\mu$ and case (A).

Also (C)  applies to $x_{\lambda}\cdot x_{\gamma}\cdot w$
since
\begin{multline*}
\eta(\lambda)+\eta(\gamma)+\eta(w)\leq
\eta(\lambda)+\eta(\gamma)+\eta(\mu)+\eta(Q)-1\\=\eta(\lambda)
+\eta(\gamma)+\eta(N)-1\leq r
\end{multline*}
and the induction hypothesis.

The preceding remarks show that (C) applies to
\[
x_{\lambda}\cdot x_{\gamma}\cdot x_{\mu}\cdot z_{Q}=
x_{\lambda}\cdot x_{\gamma}\cdot (z_{\mu}z_{Q})+x_{\lambda}\cdot
x_{\gamma}\cdot w
\]
Using \re{dblstar} and multiplying by $x_{\lambda}$,
\begin{align*}
x_{\lambda}\cdot x_{\mu}\cdot z_{N} =&x_{\lambda}\cdot T(x_{\mu}\otimes
x_{\gamma})\cdot z_{Q}+x_{\lambda}\cdot [x_{\mu},x_{\gamma}]\cdot z_{Q}\\
=&\underbrace{q_{\mu\gamma}x_{\lambda}\cdot x_{\gamma}}\cdot x_{\mu}\cdot
z_{Q}+x_{\lambda}\cdot \langle x_{\mu},x_{\gamma}\rangle \cdot z_{Q}+
x_{\lambda}\cdot [x_{\mu},x_{\gamma}]\cdot z_{Q}\\
=&q_{\mu\gamma}q_{\lambda\gamma}x_{\gamma}\cdot x_{\lambda}\cdot
x_{\mu}\cdot z_{Q}+q_{\mu\gamma}\langle  x_{\lambda},x_{\gamma}\rangle
\cdot x_{\mu}\cdot z_{Q}\\
+&q_{\mu\gamma}[x_{\lambda},x_{\gamma}]\cdot
x_{\mu}\cdot z_{Q}+x_{\lambda}\cdot\langle x_{\mu},x_{\gamma}\rangle
\cdot z_{Q}+x_{\lambda}\cdot [x_{\mu},x_{\gamma}]\cdot z_{Q}.
\end{align*}
Recall that  $\lambda$, $\mu$ are interchangeable:
\begin{align*}
x_{\mu}\cdot x_{\lambda}\cdot z_{N} &=q_{\lambda\gamma}q_{\mu\gamma}x_{\gamma}\cdot x_{\mu}\cdot
x_{\lambda}\cdot z_{Q}+q_{\lambda\gamma}\langle  x_{\mu},x_{\gamma}\rangle
\cdot x_{\lambda}\cdot z_{Q}\\
&+q_{\lambda\gamma}[x_{\mu},x_{\gamma}]\cdot
x_{\lambda}\cdot z_{Q}+x_{\mu}\cdot\langle x_{\lambda},x_{\gamma}\rangle
\cdot z_{Q}+x_{\mu}\cdot [x_{\lambda},x_{\gamma}]\cdot z_{Q}.
\end{align*}
Now use
$\eta(\lambda)+\eta(\mu)+\eta(Q)=\eta(\lambda)+\eta(N)\leq r$ to obtain
\begin{multline}
x_{\lambda}\cdot x_{\mu}\cdot z_{N}-q_{\lambda\mu}x_{\mu}\cdot
x_{\lambda}\cdot z_{N} = \\
 q_{\mu\gamma}q_{\lambda\gamma}x_{\gamma}\cdot
([x_{\lambda},x_{\mu}]+\langle x_{\lambda},x_{\mu}\rangle)\cdot
z_{Q}+\\
q_{\mu\gamma}\langle x_{\lambda},x_{\gamma}\rangle \cdot x_{\mu}\cdot
z_{Q}-q_{\lambda\mu} x_{\mu}\cdot\langle
x_{\lambda},x_{\gamma}\rangle \cdot z_{Q}+ \\
q_{\mu\gamma}[ x_{\lambda},x_{\gamma}] \cdot x_{\mu}\cdot
z_{Q}-q_{\lambda\mu} x_{\mu}\cdot[
x_{\lambda},x_{\gamma}] \cdot z_{Q} +\\
x_{\lambda}\cdot \langle x_{\mu},x_{\gamma}\rangle\cdot
z_{Q}-q_{\lambda\mu}q_{\lambda\gamma}\langle x_{\mu},x_{\gamma}\rangle\cdot
x_{\lambda}\cdot z_{Q}+\\
x_{\lambda}\cdot [ x_{\mu},x_{\gamma}]\cdot
z_{Q}-q_{\lambda\mu}q_{\lambda\gamma}[ x_{\mu},x_{\gamma}]\cdot
x_{\lambda}\cdot z_{Q}.\\  \label{redstar}
\end{multline}
Furthermore,
\begin{align}
x_{\lambda}\cdot [x_{\mu},x_{\gamma}]\cdot
& z_{Q}-q_{\lambda\mu}q_{\lambda\gamma}[x_{\mu},x_{\gamma}]\cdot
x_{\lambda}\cdot z_{Q}\nonumber\\
&=-q_{\mu\gamma}x_{\lambda}\cdot
[x_{\gamma},x_{\mu}]+q_{\lambda\mu}q_{\lambda\gamma}q_{\mu\gamma}
[x_{\gamma},x_{\mu}]\cdot x_{\lambda}\cdot z_{Q}\nonumber\\
&=q_{\mu\gamma}q_{\lambda\gamma}q_{\lambda\mu}([x_{\gamma},x_{\mu}]\cdot
x_{\lambda}\cdot z_{Q}-q_{\gamma\lambda}q_{\mu\lambda}x_{\lambda}\cdot
[x_{\gamma},x_{\mu}]\cdot z_{Q}).\label{sup}\\
\intertext{\baselineskip=22pt plus 1pt minus 1pt
If we suppose $x_{\mu}<x_{\lambda}$ then we can make use of
multiplicativity condition and since
$\eta([x_{\gamma},x_{\mu}])+\eta(\lambda)+\eta(Q)<\eta(\lambda)
+\eta(\mu)+\eta(\gamma)+\eta(Q)=\eta(\lambda)+\eta(\mu)+\eta(N)$
we obtain that \re{sup} is equal to}
&=q_{\mu\gamma}q_{\lambda\gamma}q_{\lambda\mu}[[x_{\gamma},
x_{\mu}],x_{\lambda}]\cdot z_{Q}+
q_{\mu\gamma}q_{\lambda\gamma}q_{\lambda\mu}\langle [x_{\gamma},
x_{\mu}],x_{\lambda}\rangle\cdot z_{Q}\nonumber\\
&=-q_{\lambda\gamma}q_{\lambda\mu}[[x_{\mu},x_{\gamma}],x_{\lambda}]
\cdot z_{Q}-q_{\lambda\gamma}q_{\lambda\mu}\langle [x_{\mu},x_{\gamma}],x_{\lambda}\rangle
\cdot z_{Q}. \label{red1}
\end{align}
Using multiplicativity again and since
$\eta(x_{\gamma})+\eta([x_{\lambda},x_{\mu}])+\eta(Q)<\eta(x_{\lambda})+
\eta(\mu)+\eta(\gamma)+\eta(Q)=\eta(x_{\lambda})+\eta(x_{\mu})+\eta(N)$
we may write
\begin{multline}
x_{\gamma}\cdot [x_{\lambda},x_{\mu}]\cdot
z_{Q}=\\
q_{\gamma\mu}q_{\gamma\lambda}[x_{\lambda},x_{\mu}]\cdot
x_{\gamma}\cdot z_{Q}+[x_{\gamma},[x_{\lambda},x_{\mu}]]\cdot
z_{Q}+\langle x_{\gamma},[x_{\lambda},x_{\mu}]\rangle\cdot z_{Q}
\label{red2}
\end{multline}
Sustitute \re{red1} and \re{red2} in \re{redstar},
\begin{align*}
x_{\lambda}\cdot x_{\mu}\cdot z_{N}&- q_{\mu\lambda}x_{\mu}\cdot
x_{\lambda}\cdot z_{N}=\\
 &  [x_{\lambda},x_{\mu}]\cdot x_{\gamma}\cdot z_{Q}+
 q_{\mu\gamma}q_{\lambda\gamma}[x_{\gamma},[x_{\lambda},x_{\mu}]]
 \cdot z_{Q}\\
 &+q_{\mu\gamma}q_{\lambda\gamma}\langle
 x_{\gamma},[x_{\lambda},x_{\mu}]\rangle \cdot z_{Q}
 +q_{\mu\gamma}q_{\lambda\gamma}x_{\gamma}\cdot \langle
 x_{\lambda},x_{\mu}\rangle \cdot z_{Q} \\
 &+q_{\mu\gamma}\langle x_{\lambda},x_{\gamma}\rangle \cdot x_{\mu}\cdot
 z_{Q}-q_{\lambda\mu}x_{\mu}\langle x_{\lambda},x_{\gamma}\rangle \cdot
 z_{Q}\\
 &+ q_{\mu\gamma}[ x_{\lambda},x_{\gamma}] \cdot x_{\mu}\cdot
 z_{Q}-q_{\lambda\mu}x_{\mu}[ x_{\lambda},x_{\gamma}] \cdot
 z_{Q}\\
 &+x_{\lambda}\cdot\langle x_{\mu},x_{\gamma}\rangle\cdot
 z_{Q}-q_{\lambda\mu}q_{\lambda\gamma}\langle x_{\mu},x_{\gamma}\rangle
 \cdot x_{\lambda}\cdot z_{Q}\\
 &-q_{\lambda\gamma}q_{\lambda\mu}[[x_{\mu},x_{\gamma}],x_{\lambda}]
 \cdot z_{Q}-q_{\lambda\gamma}q_{\lambda\mu}\langle
 [x_{\mu},x_{\gamma}],x_{\lambda}\rangle \cdot z_{Q}\\
\intertext{since $L$ is adequate,}
 &=[x_{\lambda},x_{\mu}]\cdot x_{\gamma}\cdot z_{Q}+\langle
 x_{\lambda},x_{\mu} \rangle \cdot x_{\gamma}\cdot z_{Q}+\\
 q_{\mu\gamma}q_{\lambda\gamma}[x_{\gamma},[x_{\lambda},x_{\mu}]&]
 \cdot z_{Q}
 -q_{\lambda\gamma}q_{\lambda\mu}[[x_{\mu},x_{\gamma}],x_{\lambda}]
 \cdot z_{Q}
 -q_{\lambda\mu}[x_{\mu},[x_{\lambda},x_{\gamma}]]\cdot z_{Q}.
\end{align*}
Thanks to Jacobi identity and (A) we get
\begin{equation}
x_{\lambda}\cdot x_{\mu}\cdot z_{N}-  q_{\lambda\mu}x_{\mu}\cdot
x_{\lambda}\cdot z_{N}-\langle x_{\lambda}, x_{\mu}\rangle \cdot
z_{N}=[x_{\lambda},x_{\mu}]\cdot z_{N}
\label{menor}
\end{equation}
if $x_{\mu}<x_{\lambda}$.

Multiplying both sides of \re{menor}   by $-q_{\lambda\mu}$ and using
antisymmetry, we get
\[
x_{\mu}\cdot x_{\lambda}\cdot z_{N}-q_{\lambda\mu}x_{\lambda}\cdot
x_{\mu}\cdot z_{N}-\langle x_{\mu},x_{\lambda}\rangle\cdot
z_{N}=[x_{\mu},x_{\lambda}]\cdot z_{N}
\]
so \re{menor} also holds  if $x_{\lambda}<x_{\mu}$.
\end{proof}
\end{section}

\begin{section}{Representations}\label{rep}
\begin{define}
Let $L$ be a basic $T$-Lie algebra and $V$ a $k$-module. A
$\,k$-morphism $\underline{\quad}\cdot \underline{\quad}:
L\otimes_{k} V\rightarrow V$ is called {\it
representation of $L$} if it satisfies
\[
x\cdot y\cdot v-T(x\otimes y)\cdot v=[x,y]\cdot v,\;\forall x,y\in
L,\,\forall v\in V
\]
where $(a\otimes b)\cdot v$ means $a\cdot b\cdot v$.
\end{define}

\begin{teo}
If $L$ is an adequate basic $T$-Lie algebra then
$L$ has a natural representation on its $q$-symmetric
algebra $\mathcal{S}(L)$.
\end{teo}

\begin{cor}
\begin{sloppypar}
If $L$ is an adequate basic $T$-Lie algebra then its
universal enveloping algebra $U(L)$ has a representation on
the $q$-symmetric algebra $\mathcal{S}(L)$.
\end{sloppypar}
\end{cor}

Inside $(sl_{n}^{+})_{q}$, $n\geq 4$, the $k$-submodules generated by the
basic elements given at Figure \ref{tipo_A_n}
\begin{figure}[h]
$$
\begin{array}{ccccc}
 e_{ij} &        & e_{jk} &      & e_{kl} \\
        & e_{ik} &        &e_{jl}         \\
        &        & e_{il}
\end{array}
$$
\caption{A basic $T$-Lie algebra of type $(sl_{4}^{+})_{q}$}\label{tipo_A_n}
\end{figure}
have a structure of basic $T$-Lie algebra that looks like
$(sl_{4}^{+})_{q}$. But such algebra have a graduation given by
$\eta(e_{ab})=a(b-a)$, this in not, in general, the graduation of
$(sl_{4}^{+})_{q}$.  However, Figure \ref{tipo_A_n} still is a basic $T$-Lie
algebra.

The basic $T$-Lie algebras given by Figure \ref{tipo_A_n} will be called {\it of type
$(sl_{4}^{+})_{q}$}. In a similar way we may define {\it basic $T$-Lie algebras of type
$(sl_{n}^{\pm})_{q}$}.
\begin{ejem}\label{sl4+}
Every basic $T$-Lie algebra of type $(sl_{4}^{\pm})_{q}$ is adequate.
\begin{proof}
Suppose $x_{\lambda}>x_{\mu}>x_{\gamma}\in\mathcal{B}$, $x_{\gamma}\leq z_{N}$
such that
$\eta(x_{\lambda})+\eta(x_{\mu})+\eta(x_{\gamma})+\eta(N)\leq r+1$. We
have to prove that \re{inducc} holds.

Note that
$\langle x_{\lambda},x_{\mu}\rangle=0$ for any
$x_{\lambda},x_{\mu}\in\mathcal{B}$ except $e_{jk},e_{il}$ , so each
term  in the equation \re{inducc} vanishes or
$e_{jk},e_{il}$ appears. This means that the equation \re{inducc} holds
trivially except in the following cases:
\[
e_{ij}<e_{jk}<e_{jl},\,e_{ij}<e_{ik}<e_{jl},\,e_{ik}<e_{jk}<e_{kl},\,e_{ik}<e_{jl}<e_{kl}
\]

\underline{Case} $e_{ij}<e_{jk}<e_{jl}$: The left side of \re{inducc} vanishes
whereas the right side is:
\begin{multline*}
e_{jk}\cdot e_{il}\cdot z_{N}-q^{2}e_{il}\cdot e_{jk}\cdot
z_{N}+q\langle e_{ik},e_{jl}\rangle\cdot z_{N} \\=e_{jk}\cdot e_{il}\cdot z_{N}
-q^{2}e_{il}\cdot e_{jk} \cdot z_{N}+(q^{2}-1)e_{il}\cdot e_{jk}\cdot
z_{N}=0,
\end{multline*}
because $e_{il}\cdot e_{jk}\cdot z_{N}=e_{jk}\cdot e_{il}\cdot z_{N}$ since
$\eta(e_{il})+\eta(e_{jk})<\eta(e_{ij})+\eta(e_{jk})+\eta(e_{jl})$ and
supposition \re{hipinducc}.

\underline{Case} $e_{ij}<e_{ik}<e_{jl}$: Let be
$d=\eta(e_{ij})+\eta(e_{ik})+\eta(e_{jl})$.
The left side of \re{inducc} is
\begin{multline*}
(q^{-1}-q)e_{il}\cdot e_{jk}\cdot e_{ij}\cdot z_{N}\\
=(q^{-1}-q)(q e_{il}\cdot e_{ij}\cdot e_{jk}\cdot z_{N}-q e_{il}\cdot
e_{ik}\cdot z_{N})
\end{multline*}
\begin{equation*}
=(q^{-1}-q)(e_{ij}\cdot e_{il}\cdot e_{jk}\cdot z_{N})-(q^{-1}-q)q e_{il}\cdot
e_{ik}\cdot z_{N}
\end{equation*}
$(\,\eta(e_{il})+\eta(e_{ij})+\eta(e_{jk})<d$ and \re{hipinducc} )
\begin{equation*}
=(q^{-1}-q)(e_{ij}\cdot e_{il}\cdot e_{jk}\cdot z_{N}-e_{il}\cdot
e_{ik}\cdot z_{N}+q e_{ik}\cdot e_{il}\cdot z_{N})
\end{equation*}
$(\,\eta(e_{ik})+\eta(e_{il})<d\,),$
and this is the right side of \re{inducc}.

The remaing cases are similar.
\end{proof}
\end{ejem}
\begin{ejem}
Every basic $T$-Lie algebra of type $(sl_{n}^{\pm})_{q}$ is adequate,
$n=5,6$.
\begin{proof}
By similar calculations as in the previous example.
\end{proof}
\end{ejem}

Note that the symbol $\cdot z_{N}$ is redundant in calculations at example
\ref{sl4+}. This remark leads to the following lemma.

Let $\otimes_{k} L$ be the tensorial $k$-algebra and $J_{r}$ the
$k$-submodule generated by
\begin{align}
x_{\alpha}\otimes x_{\beta} \otimes x_{\delta}-&T(x_{\alpha}\otimes
x_{\beta})\otimes x_{\delta}-[x_{\alpha},x_{\beta}]\otimes x_{\delta}, \\
x_{\alpha}\otimes x_{\beta}\otimes x_{\delta}-&x_{\alpha}\otimes
T(x_{\beta}\otimes x_{\delta})-x_{\alpha}\otimes [x_{\beta},x_{\delta}]
\end{align}
for $\eta(\alpha)+\eta(\beta)+\eta(\delta)\leq r,\,\forall x_{\alpha},x_{\beta},
x_{\delta}\in\mathcal{B}$.
\begin{lema}\label{adequate}
$L$ is adequate if
\begin{multline}
\langle x_{\lambda}, x_{\mu}\rangle\otimes x_{\gamma}
-q_{\lambda\mu}[x_{\mu},[x_{\lambda},x_{\gamma}]]\equiv\\
q_{\mu\gamma}q_{\lambda\gamma}\langle x_{\gamma}[x_{\lambda},x_{\mu}]
\rangle
+q_{\mu\gamma}q_{\lambda\gamma}x_{\gamma}\otimes
\langle x_{\lambda},x_{\mu}\rangle\\
+q_{\mu\gamma}\langle x_{\lambda},x_{\gamma}\rangle\otimes x_{\mu}
-q_{\lambda\mu}x_{\mu}\otimes\langle x_{\lambda},x_{\gamma}\rangle\\
+q_{\mu\gamma}[x_{\lambda},x_{\gamma}]\otimes x_{\mu}
-q_{\lambda\mu}x_{\mu}\otimes [x_{\lambda},x_{\gamma}]\\
+x_{\lambda}\otimes \langle x_{\mu},x_{\gamma}\rangle
-q_{\lambda\mu}q_{\lambda\gamma}\langle x_{\mu},x_{\gamma}\rangle\otimes
x_{\lambda}-q_{\lambda\gamma}q_{\lambda\mu}\langle
[x_{\mu},x_{\gamma}],x_{\lambda}\rangle \mod J_{r}
\label{lemainducc}
\end{multline}
for every $x_{\lambda}>x_{\mu}>x_{\gamma}\in\mathcal{B}$ such that
$\eta(\lambda)+\eta(\mu)+\eta(\gamma)\leq r+1$.
\end{lema}
\end{section}

\begin{section}{Poincar\'e-Birkhoff-Witt Theorem}\label{PBW}

Let us define
\[\mathcal{T}^{n}=\underbrace{L\otimes_{k}\ldots\otimes_{k}
L}_{n-times}
\]

For $u=x_{\lambda_{1}}\otimes\ldots\otimes
x_{\lambda_{n}}\in\mathcal{T}^{n}$
define
$\delta(u)=\eta(x_{\lambda_{1}})+\ldots+\eta(x_{\lambda_{n}})$,
$D(u)=\#\{(x_{\lambda_{i}},x_{\lambda_{j}})\,|\,
x_{\lambda_{i}}>x_{\lambda_{j}}\text{ and }i<j\}$.
If $u\in\otimes_{k}L$ and $u=\sum_{i}\xi_{i}u_{i}$ with
$u_{i}\in\mathcal{T}^{i}$, $\xi_{i}\in k,\,\forall i,$ let us put
\begin{gather}
D(u)=\max\{D(u_{i})\,|\,\xi_{i}\neq 0,\,i\}\\
\delta(u)=\max\{\delta(u_{i})\,|\,\xi_{i}\neq 0,\,i\}
\end{gather}
The number $D(u)$ is called  {\it  the disorder of $u$}.

Denote by $\mathcal{T}_{p}$  the $k$-submodule generated by
$u\in\otimes_{k}L$ such that $\delta(u)\leq p$.

\begin{define}
A sequence $(x_{\lambda_{1}},\ldots,x_{\lambda_{n}})$ of elements in a
basis of a basic $T$-Lie algebra is called {\it non-decreasing} if
$x_{\lambda_{1}}\leq\ldots\leq x_{\lambda_{n}}$ and
$x_{\lambda_{i}}=x_{\lambda_{i+1}}$ if and only if
$S(x_{\lambda_{i}}\otimes x_{\lambda_{i+1}})=
x_{\lambda_{i}}\otimes x_{\lambda_{i+1}}$.
\end{define}

\begin{teo}[Poincar\'e-Birkhoff-Witt]
Let $L$ be an adequate $T$-Lie algebra with a basis $\mathcal{B}$. The
monomials  formed by finite non-decreasing sequences of elements in
$\mathcal{B}$ constitute a free $k$-basis of the  universal enveloping algebra
$U(L)$.
\end{teo}
\begin{proof}

Let $P:\otimes_{k}L\rightarrow U(L)$ be the canonical $k$-morphism,
$\mathcal{M}$ the $k$-submodule generated by the monomials described
in the formulation of the theorem.
We have to prove that
$U(L)=\mathcal{M}$.  Note that
\[
U(L)=\sum_{p=1}^{\infty}P(\mathcal{T}_{p})
\]
If $p=1$ then $\mathcal{T}_{p}\subseteq L$ it follows
$P(\mathcal{T}_{p})\subseteq \mathcal{M}$. Suppose
$P(\mathcal{T}_{r})\subseteq\mathcal{M}$. It suffices to show that
$P(\mathcal{T}_{r+1})\subseteq\mathcal{M}$.

Define $\mathcal{T}^{u}_{r}$ as the $k$-submodule of
$\mathcal{T}_{r}$ generated by elements with disorder $\leq u$,
and proceed by a second induction
on the disorder. We have $P(\mathcal{T}_{r+1}^{0})\subset\mathcal{M}$.
Suppose $v=a\otimes x\otimes y\otimes b\in\mathcal{T}^{u}_{r+1}$ where
$x>y\in\mathcal{B},$ and $a\in\mathcal{T}^{n},\,b\in\mathcal{T}^{m}$
monomials form by basic elements in $\mathcal{B}$.
Then
\begin{align*}
P(v)&=P(a\otimes q_{xy}y\otimes x\otimes b)+P(a\otimes\langle
x,y\rangle\otimes b)+P(a\otimes [x,y]\otimes b)\\
    &\equiv P(a\otimes q_{xy}y\otimes x\otimes b)\mod{\mathcal{T}_{r}}
\end{align*}
but $P(a\otimes q_{xy}y\otimes x\otimes b)\in
\mathcal{T}^{u-1}_{r+1}\subseteq \mathcal{M}$. Hence
$P(v)\in\mathcal{M}$ it follows $P(\mathcal{T}_{r+1})\subseteq \mathcal{M}$.

It remains to prove linear independence. For a given sequence
$\Sigma=(x_{\lambda_{1}},\ldots,x_{\lambda_{n}})$ of
non-decreasing elements of $\mathcal{B}$, define
$x_{\Sigma}=x_{\lambda_{1}}\dots\, x_{\lambda_{n}}\newline\in U(L)$.
Suppose
\[
\sum_{i}\xi_{i}x_{\Sigma_{i}}=0
\]
where each $\Sigma_{i}$ is a sequence non-decreasing and
$\xi_{i}\in k,\,\forall i$.
Using the representation of $U(L)$, we get from lemma (C)
\[
0= \sum_{i}\xi_{i}x_{\Sigma_{i}}\cdot 1
 =\sum_{i}\xi_{i}z_{\Sigma_{i}}
\]
and because the linear independence of the
$z_{\Sigma_{i}}\in\mathcal{S}(L)$, it
follows that $\xi_{i}=0,\,\forall i$.
\end{proof}
\begin{cor}
\begin{sloppypar}
$U(sl_{n}^{\pm})_{q}$ has a basis of the Poincar\'e-Birkhoff-Witt type,
$n=2,3,4,5$.
\end{sloppypar}
\end{cor}
\end{section}

\begin{section}{Braids}\label{braids}

\begin{prop}
In $L=(sl_{n}^{\pm})_{q}$, $n=2,3,4$  it holds the braid equation:
\[
T_{12}T_{23}T_{12}|_{\,^{3}L}=T_{23}T_{12}T_{23}|_{\,^{3}L},
\]
where
$T_{12}=T\otimes_{k}Id_{L},T_{23}=Id_{L}\otimes_{k}T$.
\end{prop}
\begin{proof}
By straightforward calculations on the basic elements.
(\,Using {\it Mathematica} \Wolfram).
\end{proof}
\begin{prop}
The presymmetry $S$ of a $T$-Lie algebras holds the braid equation
\[
S_{12}S_{23}S_{12}=S_{23}S_{12}S_{23}
\]
\end{prop}
\begin{proof}
Let $x,y,z$ be basic elements. Then
\[
S_{12}S_{23}S_{12}(x\otimes y\otimes
z)=q_{x,y}q_{x,z}q_{y,z}z\otimes y\otimes x
= S_{23}S_{12}S_{23}(x\otimes y\otimes z)
\]
\end{proof}

\begin{nota}
\begin{sloppypar}
The symmetry of $\widetilde{(sl_{4}^{+})}_{q}$  is a braid morphism, however
we have no PBW theorem
for $U\widetilde{(sl_{4}^{+})}_{q}$. As a
consecuence the PBW theorem is independent from the
braid equation.
\end{sloppypar}
\end{nota}

\end{section}

\begin{section}{Non-Standard Quantum Deformations of
$GL(n)$}\label{nonstandard}
\begin{define}
Let $p,q$ be units in a commutative unitary ring $k$  with $pq\neq 1$ and
choose $\alpha(\alpha-1)/2$ discrete parameters $\epsilon_{ij}$,
$\epsilon_{ij}=\pm 1,\,1\leq i<j\leq
\alpha,\,\epsilon_{ii}=1,\,\epsilon_{ji}=\epsilon_{ij}$. Let $m,n$ be
positive
integers such that  $m,n\leq \alpha$

The $k$-module $L_{p,q,\epsilon}(n,m,k)$ is then defined to be the free
$k$-module with basis
\[
\mathcal{B}=\{Z_{i}^{j}\,|\,1\leq i\leq n,\,1\leq j\leq m\}.
\]
Now define an order on $\mathcal{B}$ and morphisms
$S,\,T\,,\langle,\rangle,\,[,]$ copying  the structure of $L(n,k)$ in
example \ref{nonstand}.
\end{define}

\begin{prop}
$L_{p,q,\epsilon}(n,m,k)$ has a structure of basic $T$-Lie algebra.
\end{prop}

In a similar way to the algebras of type $(sl_{n}^{+})_q$, (see section \ref{rep})
we can define {\it algebras of type $L_{p,q,\epsilon}(n,m,k)$.}
\begin{lema}
Every algebra of type $L_{p,q,\epsilon}(\lambda,\mu,k)$ is an adequate basic $T$-Lie algebra,
where $\lambda,\mu\in\{2,3\}$.
\end{lema}

\newcommand{\Z}[2]{Z_{#1}^{#2}}

\begin{lema}
If $\Z{u}{v}>\Z{i}{j}>\Z{a}{b}$ then there exists $L$ being a $T$-Lie subalgebra
of $L_{p,q,\epsilon}(n,m,k)$  and numbers $\lambda,\mu\in\{2,3\}$
such that $L$ is of type $L_{p,q,\epsilon^{\prime}}(\lambda,\mu,k)$ and
$\{\,\Z{u}{v},\Z{i}{j},\Z{a}{b}\,\}\subset L$.
\end{lema}
\begin{proof}

Let us put the basic elements in a matrix array (Figure \ref{d} (a)).
\begin{figure}[h]
$$
\text{(a)}
\begin{array}{cccc}
\Z{1}{1}& \Z{1}{2}&\ldots &\Z{1}{n} \\
\vdots  & \vdots  &       &\vdots   \\
\Z{m}{1}& \Z{m}{2}&\ldots &\Z{m}{n}
\end{array}
\quad
\text{(b)}
\begin{array}{cc}
\circ_{i}^{j} & \circ_{i}^{j+u}\\
\circ_{i+u}^{j} & \circ_{i+u}^{j+u}
\end{array}
$$
\caption{\baselineskip=22pt plus 1pt minus 1pt
(a) The basic $T$-Lie algebra
$L_{p,q,\epsilon^{\prime}}(\lambda,\mu,k)$. (b) Diagonal relationship.}
\label{d}
\end{figure}

Note that for  positive integers $u,v$ the elements appearing in the
pseudobracket definition are in a diagonal relationship (Figure
\ref{d}(b)\,),
and they form a free basis of a $T$-Lie algebra of type
$L_{p,q,\epsilon^{\prime}}(2,2,k)$, where
$\epsilon^{\prime}=\{1,\epsilon_{i,i+u},\,\epsilon_{j,j+u}\}$.

For $\Z{u}{v}>\Z{i}{j}>\Z{a}{b}$ there are several cases. The cases
given by  Figure \ref{ult}(a), \ref{ult}(b), \ref{ult}(c),
\begin{figure}[t]
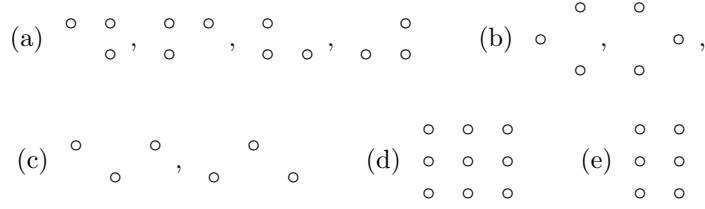

$$
\text{(a)}\;
\begin{array}{cc}
\circ & \circ \\
      & \circ
\end{array},
\;
\begin{array}{cc}
\circ & \circ \\
\circ &
\end{array},
\;
\begin{array}{cc}
\circ \\
\circ & \circ
\end{array},
\;
\begin{array}{cc}
      & \circ \\
\circ & \circ
\end{array}
\qquad
\text{(b)}\;
\begin{array}{cc}
      &\circ \\
\circ & \\
      &\circ
\end{array},\;
\begin{array}{cc}
\circ & \\
      &\circ \\
\circ
\end{array},
$$

$$
\text{(c)}\;
\begin{array}{ccc}
\circ & &\circ\\
      &\circ
\end{array},\;
\begin{array}{ccc}
     &\circ &\\
\circ&      &\circ
\end{array}
\qquad
\text{(d)}\;
\begin{array}{ccc}
\circ & \circ & \circ \\
\circ & \circ &\circ \\
\circ & \circ &\circ \\
\end{array}
\qquad
\text{(e)}\;
\begin{array}{cc}
\circ &\circ \\
\circ &\circ \\
\circ &\circ
\end{array}
$$
\caption{Some cases in $L_{p,q,\epsilon}(n,m,k)$}\label{ult}
\end{figure}
or they form a triangle which can be fitted, with vertices on the border,
inside of the rectangle at Figure \ref{ult}(d).

In the case given by Figure \ref{ult}(a) we may complete each triangle to a square and
obtain $L_{p,q,\epsilon_{0}}(2,2,k)$. In the case given by Figure
\ref{ult}(b), each triangle
can be completed to a rectangle in the form of Figure \ref{ult}(e)
and we get $L_{p,q,\epsilon_{1}}(3,2,k)$. Similarly in the case given by
Figure \ref{ult}(c) we get $\,L_{p,q,\epsilon_{2}}(2,3,k)$. Finally, in the
case given by  Figure
\ref{ult}(d), we obtain $L_{p,q,\epsilon_{3}}(3,3,k)$.
\end{proof}

\begin{teo}
$L_{p,q,\epsilon}(n,m,k)$   is an adequate basic $T$-Lie algebra.
\end{teo}

\begin{cor}
The monomials formed by non-decreasing finite sequences of elements in
\[
\mathcal{B}=\{Z_{i}^{j}\,|\,1\leq i\leq n,\,1\leq j\leq m\}
\]
constitute a free basis of the $k$-module $M_{p,q,\epsilon}(n,m,k)=
U\,L_{p,q,\epsilon}(n,m,k)$.
\end{cor}
\end{section}

\begin{section}{$U(sl_{n+1}^{+})_{q}$. The general case.}\label{A+}
\begin{lema}\label{lemJ}
Let $e_{ab},e_{uv},e_{ij}$ be basic elements in $(sl_{n+1}^{+})_{q}$ and
$[,]$ usual bracket in $sl_{n+1}$.
\begin{enumerate}
\item $e_{ab}<e_{uv}$ if and only if $a+b<i+j$ or $a+b=i+j$ and $b<j$.
\item If $S(e_{ab}\otimes e_{uv})=q^{c_{ab,uv}}e_{uv}\otimes e_{ab}$ and
$e_{ab}<e_{uv}$ then
\[
c_{ab,uv}=-\delta_{v,a}+\delta_{v,b}+\delta_{u,a}-\delta_{u,b}
\]
\item  If  $e_{ab}<e_{uv}<e_{ij}$ then
\[
q^{c_{uv,ab}}[\,e_{uv},[e_{ab},e_{ij}]\,]=[\,e_{uv},[e_{ab},e_{ij}]\,]_{q}
\]
\end{enumerate}
\end{lema}
\begin{proof}
\begin{enumerate}
\item By the order definition.
\item It follows from the formula
$[e_{ij},e_{kl}]=\delta_{jk}e_{il}-\delta_{li}e_{kj}$ in the classical
Lie algebra $sl_{n}$.
\item \begin{sloppypar}Since $e_{ab}<e_{ij}$ we can suppose $b=i$. We have to prove
$q^{c_{uv,ab}}[e_{uv},e_{aj}]=[e_{uv},e_{aj}]_{q}$. There are two cases
     \end{sloppypar}
     \begin{enumerate}
     \item\label{cas1} $e_{uv}<e_{aj}$;
     \item\label{cas2} $e_{uv}>e_{aj}$.
     \end{enumerate}

\eqref{cas1}: If $[e_{uv},e_{aj}]\neq 0$ then $v=a$ and $u<v=a<b$, it
follows $e_{uv}<e_{ab}$ since $u+v<a+b$. A contradiction. Therefore
$[e_{uv},e_{aj}]_{q}=[e_{uv},e_{aj}]=0$. 

\eqref{cas2}: We have to prove
\[
q_{uv,ab}[e_{aj},e_{uv}]=q_{uv,aj}[e_{aj},e_{uv}]
\]
Both sides are zero because if not then $j=u$ and $i<j=u<v$, these imply
$i+j<v+u$, and then $e_{ij}<e_{uv}$. Again, we have a contradiction.
\end{enumerate}
\end{proof}
\begin{teo}
$(sl_{n+1}^{+})_{q}$ is an adequate basic $T$-Lie algebra.
\end{teo}
\begin{proof}
Let $\mathcal{B}$ be the canonical basis of $sl_{n+1}$, and write the basic elements
of $\mathcal{B}$ in the form $e_{ij}$. Now, we put this basic elements
in an upper triangular array (Figure \ref{diagrama}(a)\,).
Note that, if $\langle e_{ij},e_{(i+u)(j+v)}\rangle\neq 0$ then the
elements appearing in the pseudobracket definition are in a diagonal
relationship (Figure \ref{diagrama}(b)\,),
and if $[e_{ij},e_{(i+u)(j+v)}]_{q}\neq 0$ then $j=i+u$ and we get  the
Figure \ref{diagrama}(c).

So, if we suppose $e_{ij}>e_{uv}>e_{ab}$ then the elements appearing
in the formulation of lemma \ref{adequate} (brackets and pseudobrackets)
can be fitted inside of a square of the form of Figure \ref{diagrama}(d),
and such square can be extended to an upper triangle (Figure
\ref{diagrama}(e)),
but this triangle gives a strictly graded algebra of type  $(sl_{6}^{+})_{q}$.
 Since these algebras  satisfies the condition of lemma \ref{adequate},
in particular the elements $e_{ij}>e_{uv}>e_{ab}$ satisfies this
condition.  Besides,
\begin{multline*}
[,]_{q}(\:(Id\otimes[,]_{q})S_{12}S_{23}-([,]_{q}\otimes
Id)S_{23}S_{12}+\\
(Id\otimes [,]_q)S_{23}S_{12}\,)(e_{ij}\otimes e_{uv}\otimes
e_{ab})=
\end{multline*}
\begin{gather*}
q^{c_{uv,ab}+c_{ij,ab}+c_{ij,uv}}([[e_{ab},e_{uv}],e_{ij}]
-[e_{ab},[e_{uv},e_{ij}]]+[e_{uv},[e_{ab},e_{ij}]])=0
\end{gather*}
since lemma \ref{lemJ} and the Jacobi identity in $sl_{n}^{+}$.

We conclude that $(sl_{n+1}^{+})_{q}$ is an adequate basic $T$-Lie algebra.
\end{proof}

\begin{figure}[t]
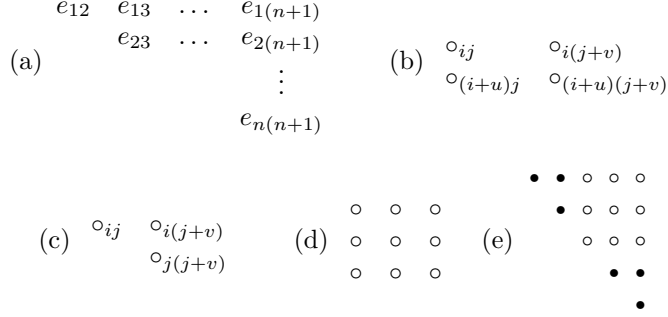

$$
\text{(a)}
\begin{array}{cccc}
e_{12}&e_{13}&\ldots&e_{1(n+1)}\\
      &e_{23}&\ldots&e_{2(n+1)}\\
      &      &      &\vdots\\
      &      &      &e_{n(n+1)}
\end{array}
\qquad
\text{(b)}\;
\begin{array}{lll}
\circ_{ij}     & \circ_{i(j+v)} \\
\circ_{(i+u)j} & \circ_{(i+u)(j+v)}
\end{array}
$$

$$
\text{(c)}\;
\begin{array}{lll}
\circ_{ij}     & \circ_{i(j+v)} \\
               & \circ_{j(j+v)}
\end{array}
\qquad
\text{(d)}\;
\begin{array}{cccc}
\circ&\circ&\circ\\
\circ&\circ&\circ\\
\circ&\circ&\circ
\end{array}
\newcommand{\cir}{\put(0,2.5){\circle*{3.2}}}
\newcommand{\ci}{\put(0,2.5){\circle{3.2}}}
\quad
\text{(e)}\;
\begin{array}{cccccc}
\cir&\cir&\ci&\ci&\ci\\
    &\cir&\ci&\ci&\ci\\
    &    &\ci&\ci&\ci\\
    &    &   &\cir&\cir\\
    &    &   &   &\cir
\end{array}
$$
\caption{Some cases in $(sl_{n+1}^{+})_{q}$}\label{diagrama}
\end{figure}

\begin{lema}\label{lema-An}
Suppose $e_{ij}<e_{ab}\in U_{q}^{+}(sl_{n+1})$. Then the following
equations are satisfied in $U_{q}^{+}(sl_{n+1})$,
\begin{equation}
[e_{ij},e_{ab}]=\begin{cases}
               e_{ij}e_{ab}-q e_{ab}e_{ij},\text{ if }i=a\text{ or
               }j=b\\
               e_{ij}e_{ab}-e_{ab}e_{ij}-\langle
               e_{ij},e_{ab}\rangle\text{ if }i\neq a,
               j\neq b\text{ and }j\neq a, \\
               e_{ij}e_{ab}-q^{-1}e_{ab}e_{ij},\text{ if }j=a.
               \end{cases}\label{ec}
\end{equation}
\end{lema}
\begin{proof}
By induction on $n$.
For the cases $n=1,2,3,4,5$ the equations \ref{ec} can be verified by
straightforward calculations.
So we may suppose $n>5$. Let us consider the Figure \ref{diagrama}(a).
Such diagram can be thought as formed by two  overlaping triangles. The
first one, a triangle $T_{1}$ with vertices $e_{12},e_{1n},e_{(n-1)n}$ and
the second one, a triangle $T_{2}$ with vertices
$e_{23},e_{2(n+1)},e_{n(n+1)}$.

The elements in $T_{i}$ generate a $k$-subalgebra  isomorphic to
$U_{q}^{+}(sl_{n})$, $i=1,2$. Then, if $e_{ij}$ and $e_{ab}$ are both in
$T_{1}$ or $T_{2}$, the equations \re{ec} holds. As a consecuence, we
may suppose $i=1$ and $b=n+1$, and put $j\neq n+1$ and $a\neq 1$.

At the Figure \ref{diagrama}(a) join the node $rs$ with the node $uv$ if
$[e_{rs},e_{uv}]_{q}\neq 0$. We have several cases given by Figure
\ref{cases}, (in the first and third cases, since $e_{12},e_{2j},e_{nn}$ are in $T_{1}$
and the induction
hypothesis there is not arrow between $12$ and $an$, whereas
there is not arrow between $2j$ and $n(n+1)$ because
$e_{2j},e_{an},e_{n(n+1)}$ are in $T_{2}$).

At the first case we get a graph of type $A_{4}$, then
$e_{1j}=[e_{12},e_{2j}]_{q},$ $e_{a(n+1)}=[e_{an},e_{n(n+1)}]_{q}$
are in a subalgebra isomorphic to
$U_{q}^{+}(sl_{5})$, it follows,
\[
[e_{1j},e_{j(n+1)}]_{q}=e_{1j}e_{j(n+1)}-q^{-1}e_{j(n+1)}e_{ij}
\]
In the second case we get a graph of type $A_{3}$ then
$e_{1j}=[e_{12},e_{2n}]_{q},$ $e_{2(n+1)}=[e_{2n},e_{n(n+1)}]_{q}$ are in a
subalgebra isomorphic to $U_{q}^{+}(sl_{4})$, besides
\[
[e_{1j},e_{a(n+1)}]_{q}=e_{1n}e_{2(n+1)}-e_{2(n+1)}e_{1n}
-(q-q^{-1})e_{1(n+1)}e_{2n}
\]

\begin{figure}[t]
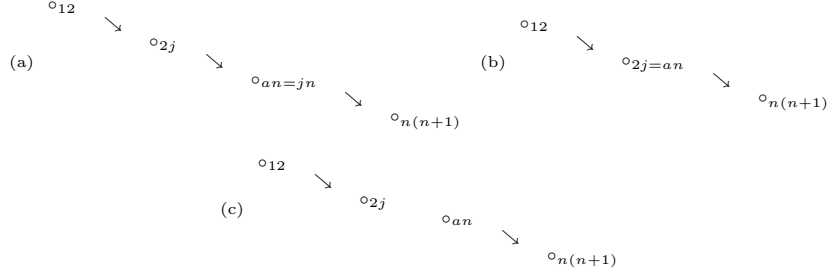

{\tiny
\[
\text{(a)}
\begin{array}{ccccccc}
\circ_{12}&      &\\
          &\searrow&\\
          &        &\circ_{2j}&\\
          &        &          &\searrow\\
          &        &          &         &\circ_{an=jn}&\\
          &        &          &         &          &\searrow\\
          &        &          &         &          &         &\circ_{n(n+1)}
\label{a}
\end{array}
\;\text{(b)}
\begin{array}{ccccc}
\circ_{12}& \\
          &\searrow \\
          &         &\circ_{2j=an}&\\
          &         &             &\searrow\\
          &         &             &        &\circ_{n(n+1)}
\label{b}
\end{array}
\]
\[
\text{(c)}
\begin{array}{ccccccc}
\circ_{12}&      &\\
          &\searrow&\\
          &        &\circ_{2j}&\\
          &        &          &         &\circ_{an}&\\
          &        &          &         &          &\searrow\\
          &        &          &         &          &         &\circ_{n(n+1)}
\label{c}
\end{array}
\]
}
\caption{(a) First case (b) second case (c) third case.}\label{cases}
\end{figure}

In the third case we may insert the node $ja$ in order to obtain
\[
\circ_{12}\rightarrow\circ_{2j}\rightarrow\circ_{ja}\rightarrow\circ_{an}
\rightarrow\circ_{n(n+1)}
\]
this graph is of type $A_{5}$, then
$e_{1j}=[e_{12},e_{2j}]_{q},e_{a(n+1)}=[e_{an},e_{n(n+1)}]_{q}$ are in a
subalgebra isomorphic to $U_{q}(sl_{6}^{+})$ and
\[
0=[e_{1j},e_{a(n+1)}]_{q}=e_{1j}e_{a(n+1)}-e_{a(n+1)}e_{1j}
\]

Now only remains the cases $e_{1j}=e_{1(n+1)}, e_{a(n+1)}=
e_{1(n+1)}$. Suppose $e_{1j}=e_{1(n+1)}$. Since
$e_{12}e_{a(n+1)}=e_{a(n+1)}$,
$e_{2(n+1)}e_{a(n+1)}=qe_{a(n+1)}e_{2(n+1)}$ and
$e_{1(n+1)}=e_{12}e_{2(n+1)}-q^{-1}e_{2(n+1)}e_{12}$ it follows,
\[
e_{1(n+1)}e_{a(n+1)}=qe_{a(n+1)}e_{1(n+1)}.
\]
In a similar way, if $e_{a(n+1)}=e_{1(n+1)}$, we get
\[
e_{1j}e_{a(n+1)}=qe_{a(n+1)}e_{1j}.
\]
\end{proof}
\begin{teo}There exists an isomorphism
\[
U_{q}^{+}(sl_{n+1})\simeq U(sl_{n+1}^{+})_{q}
\]
of $k$-algebras.
\end{teo}
\begin{proof}
Let us put $c_{ab,cd}=c_{uv}$ where $x_{u}=e_{ab},x_{v}=e_{cd}$, and
$x_{u}<x_{v}$,  $1\leq a,b,c,d\leq n+1$. From
\[
[e_{ij},e_{kl}]=\delta_{jk}e_{il}-\delta_{li}e_{kj}
\]
it follows, if $e_{ab}<e_{cd}$,
\[
c_{ab,cd}=\begin{cases}
          1,\text{ if }a=c\text{ or }b=d,\\
          0,\text{ if }a\neq c\text{ and }b\neq d,\\
          -1,\text{ if }b=c,
          \end{cases}
\]

Now use lemma \ref{lema-An} in order to obtain the following equations
in $U_{q}^{+}(sl_{n+1})$,
\[
e_{ab}e_{cd}-q^{c_{ab,cd}}e_{cd}e_{ab}=[e_{ab},e_{cd}]_{q}+\langle
e_{ab},e_{cd}\rangle,
\]
for all $1\leq a,b,c,d \leq n+1$.

We conclude $U_{q}^{+}(sl_{n+1})\simeq U(sl_{n+1}^{+})_{q}$.
\end{proof}

\begin{cor}\label{defor}
\begin{enumerate}
\item The monomials formed by non-decreasing finite sequences of elements in
\[
\mathcal{B}=\{e_{ij}\,|\,1\leq i<j\leq m\}
\]
constitute a free basis of the $k$-module $U_{q}^{+}(sl_{n+1})$, where
$m=n(n+1)/2$.
\item We have
\[(sl_{n+1}^{+})_{q}|_{t=0}=sl_{n+1}^{+}\]
and $(sl_{n+1}^{+})_{q}$  is a
deformation of $sl_{n+1}^{+}$ in the category of $T$-Lie algebras.
\end{enumerate}
\end{cor}
\end{section}

\bigskip
{\Large {\bf Acknowledgments}}
\nopagebreak
\bigskip
\nopagebreak

I would like to thank M. Durdevi\'c, R. Bautista  for their
encouragement, to C. Ramirez for giving me access to his papers collection,
A. Sudbery for drawing my attention to his own work about quantum
algebras, and thanks to J. Garc\'{\i}a-Fern\'andez, P. Garc\'{\i}a and L.
Montgomery for the computer freeware.

The computer programs used in this
work were developed by the author in collaboration with  Yaoc\'{\i}huatl E.
Arroyo.

\end{document}